\documentstyle[pra,aps,eqsecnum,amssymb]{revtex}

\draft

\begin{document}

\title{Applications of commuting difference operators
to \\ orthogonal polynomials in several variables\thanks{Work
supported by the Japan Society for the Promotion of Science
(JSPS) and by a Monbusho Grant-in-Aid.}}

\author{J. F. van Diejen}
\address{Department of Mathematical Sciences, University of Tokyo,\\
Komaba~3-8-1, Meguro-ku, Tokyo~153, Japan}

\date{December 1995}

\maketitle

\begin{abstract}
Elementary properties of the Koornwinder-Macdonald multivariable
Askey-Wilson polynomials are discussed.
Studied are the orthogonality, the difference equations, the recurrence
relations, and the orthonormalization constants for these polynomials.
Essential in our approach are certain commuting
difference operators simultaneously diagonalized by the polynomials.
\end{abstract}

\vspace{2ex}
\begin{center}
{\em Mathematics Subject Classifications (1991).} 33D45, 33D80
\end{center}

\section{Introduction}\label{sec1}
It is well-known that large part of the classical
theory of hypergeometric orthogonal
polynomials is intimately connected with the representation theory of
(simple, rank one) Lie groups and that, similarly, {\em basic}
hypergeometric orthogonal polynomials are connected with (the representation
theory of) quantum groups \cite{vil-kli:representations1}.
Moreover, it turns out that these relations with representation
theory provide a fruitful framework for generalizing
such polynomials to several
variables \cite{vil-kli:representations2}.

Around ten years ago, Askey and Wilson introduced a now famous family
of basic hypergeometric orthogonal
polynomials in one variable, which contains
many of the (basic) hypergeometric orthogonal polynomials studied in the
literature as special
(limiting) cases \cite{ask-wil:some,koe-swa:askey-scheme}.
Recently, a multivariable generalization of the Askey-Wilson
polynomials was found, first for special parameters
by Macdonald \cite{mac:orthogonal} and then
by Koornwinder \cite{koo:askey} for general parameters.
It turns out that from the viewpoint of
representation theory (and for special values of the
parameters) the Koornwinder-Macdonald multivariable
Askey-Wilson polynomials correspond to zonal spherical functions on quantum
symmetric spaces of classical type \cite{nou-sug:quantum}.
Furthermore, the polynomials may be interpreted physically
as the eigenfunctions for a
Ruijsenaars type difference version of
the Calogero-Sutherland $n$-particle model
related to the root system $BC_n$ \cite{die:diagonalization,rui:complete}.

Here we will study some elementary properties
of these
multivariable Askey-Wilson polynomials, with the
aim of
generalizing part of the well-known results by Askey and Wilson to the
case of several variables.
Attention will be mainly focussed on describing the difference equations,
the recurrence relations, and the orthonormalization constants for the
polynomials. An essential role in our discussion is played by
a family of previously introduced commuting
difference operators that are simultaneously diagonalized by the polynomials.
Most results are stated without complete proof.
A more detailed treatment including proofs can be found in
\cite{die:self-dual}.

\section{Multivariable Askey-Wilson polynomials}\label{sec2}
A natural basis for the algebra of permutation invariant and even
trigonometric polynomials is given by the monomial symmetric functions
\begin{equation}
m_\lambda (x) \; =\sum_{\lambda^\prime \in W \lambda}
e^{\alpha \sum_{j=1}^n \lambda_j^\prime x_j},
\;\;\;\;\;\; \lambda \in \Lambda = \{ \lambda \in {\Bbb Z}^n \; |\; \lambda_1
\geq \lambda_2 \geq \cdots \geq \lambda_n \geq 0\; \} , \label{monomial}
\end{equation}
where the summation is over the orbit of
$\lambda$ under the action of
the (Weyl) group $W$($\cong S_n \ltimes ({\Bbb Z}_2)^n$)
generated by permutations and sign flips
of the vector components $\lambda_1,\ldots ,\lambda_n$.
One can (partially) order the basis of monomial symmetric functions
by defining for all $\mu, \lambda \in \Lambda$
\begin{equation}\label{ord}
\mu \leq \lambda \;\;\;\;\;\;\;\; \text{iff} \;\;\;\;\;\;\;\;
\sum_{1\leq j\leq m}\mu_j \leq \sum_{1\leq j\leq m} \lambda_j
\;\;\;\;\; \text{for} \;\;\;\;\; m=1,\ldots ,n
\end{equation}
(and $\mu < \lambda$ iff $\mu \leq \lambda$ and
$\mu \neq \lambda$).

The Koornwinder-Macdonald multivariable Askey-Wilson
polynomials $p_\lambda(x)$ (with
$\lambda\in\Lambda$) are now defined as the (unique) trigonometric
polynomials satisfying
\begin{itemize}
\item[i.] $\displaystyle p_\lambda (x) = m_\lambda (x) +
\sum_{\mu \in \Lambda, \mu < \lambda }\;
 c_{\lambda ,\mu }\: m_{\mu}(x)$,
\ \ \ \ \ $\displaystyle c_{\lambda ,\mu}\in {\Bbb C}$;
\item[ii.]
$\displaystyle  \langle p_\lambda  ,  m_{\mu}\rangle_\Delta  =0$
\ \ if \ \  $\displaystyle \mu < \lambda$.
\end{itemize}
Here $\langle\cdot ,\cdot\rangle_\Delta$ denotes the $L^2$ inner product
\begin{equation}\label{ip}
\langle m_\lambda ,m_{\mu} \rangle_\Delta =
\left( \frac{\alpha}{2\pi} \right)^n
\int_{-\pi/\alpha}^{\pi/\alpha}\!\!\!\!\!\!\cdots
               \int_{-\pi/\alpha}^{\pi/\alpha}
m_\lambda (ix)\, \overline{m_{\mu}(ix)}\,  \Delta (ix)\,
dx_1\cdots dx_n
\end{equation}
with weight function
\begin{equation}\label{weight}
\Delta (x)\; =\; \Delta^+ (x)  \Delta^+ (-x),
\end{equation}
where
\begin{eqnarray}\label{weightp}
 \Delta^+ (x) &=&\prod_{1\leq j<k \leq n }
d^+_v (x_j+x_{k})\, d^+_v (x_j-x_{k})\;
\prod_{1\leq j\leq n}
d^+_w (x_j),\\
d^+_v(z)& =&q^{-gz/2}
\frac{(q^{z}; q)_\infty}
     {( q^g q^{z}; q)_\infty},\;\;\;\;\;\;\;\;\;\;\;\;\;\;
     q=e^{-\alpha },  \nonumber  \\[1ex]
   d^+_w(z)&=& q^{-(g_0+\cdots +g_3)z/2}
\frac{(q^{z},\,
      -q^{z},\,
       q^{1/2} q^{z},\,
      -q^{1/2} q^{z};\,
                                    q)_\infty    }
     {(q^{g_0} q^{z},\,
      -q^{g_1} q^{z},\,
       q^{(g_2+1/2)}q^{z},\,
      -q^{(g_3+1/2)} q^{z};\,
                                    q)_\infty    }  .\nonumber
\end{eqnarray}
(As usual $(a;q)_\infty\equiv\prod_{l=0}^\infty (1-a q^l)$ and
$(a_1,\ldots ,a_k;q)_\infty \equiv (a_1;q)_\infty \cdots (a_k;q)_\infty$.)
To ensure the convergence of the infinite products in $\Delta $
(\ref{weight}) it will be assumed that $\alpha $ is positive (so $0<q<1$); in
addition, we will also assume that $g,g_r \geq 0$, $r=0,1,2,3$.

\section{Difference equations}\label{sec3}
In \cite{die:commuting} (see also \cite[Section 7.1]{die:self-dual})
it was shown that the polynomials $p_\lambda (x)$ satisfy a system
of difference equations having the structure of eigenvalue equations
\begin{equation}\label{diffeqr}
(D_r\, p_\lambda )(x) = E_r (\rho +\lambda )\:  p_\lambda (x)
\end{equation}
for $n$ independent commuting difference operators $D_1,\ldots ,D_n$
of the form
\begin{equation}\label{ados}
D_r=\sum_{\stackrel{ J\subset \{ 1,\ldots ,n\} ,\, 0\leq |J|\leq r }
         {\varepsilon_j=\pm 1,\; j\in J}}   \!\!\!\!
U_{J^c,\, r-|J|}(x)\,  V_{\varepsilon J,\, J^c}(x)\,
T_{\varepsilon J },\;\;\;\;\;\;\;\;\;
r=1,\ldots ,n,
\end{equation}
with $(T_{\varepsilon J }f)(x)=f(x+e_{\varepsilon J })$ where
$e_{\varepsilon J }=\sum_{j\in J} \varepsilon_j e_j$.
(Here $e_j$ denotes the
$j$th unit vector in the standard basis of ${\Bbb R}^n$.)
The coefficients of the difference operators read explicitly
\begin{eqnarray*}
V_{\varepsilon J,\, K}(x) &=&
\prod_{j\in J} w(\varepsilon_jx_j) 
\prod_{\stackrel{ j,j^\prime \in J}{j<j^\prime}}
v(\varepsilon_jx_j+\varepsilon_{j^\prime}x_{j^\prime})
 v(\varepsilon_jx_j+\varepsilon_{j^\prime}x_{j^\prime}+1 )\\
& & \times  \prod_{\stackrel{j\in J}{k\in K}} v(\varepsilon_j x_j+x_k)
v(\varepsilon_j x_j -x_k),  \\ [1ex]
U_{K,p}(x)&=& (-1)^p \!\! \sum_{\stackrel{L\subset K,\, |L|=p}
                               {\varepsilon_l =\pm 1,\; l\in L}} \;
\prod_{l\in L} w(\varepsilon_l x_l)\,
\prod_{\stackrel{l,l^\prime \in L}{l<l^\prime}}
v(\varepsilon_lx_l+\varepsilon_{l^\prime}x_{l^\prime})
v(-\varepsilon_lx_l-\varepsilon_{l^\prime}x_{l^\prime}-1 ) \\
& & \times \prod_{\stackrel{l\in L}{k\in K\setminus L}}
v(\varepsilon_l x_l+x_k) v(\varepsilon_l x_l -x_k) ,
\end{eqnarray*}
with
\begin{eqnarray*}
v(z)&=& \frac{\text{sh} \frac{\alpha}{2} (g + z)}
                 {\text{sh} (\frac{\alpha}{2} z)},    \\
w(z)&=& \frac{\text{sh} \frac{\alpha}{2} (g_0 + z)}
           {\text{sh} (\frac{\alpha}{2} z)}
       \frac{\text{ch} \frac{\alpha}{2} (g_1 + z)}
            {\text{ch} (\frac{\alpha}{2} z)}
      \frac{\text{sh} \frac{\alpha}{2} (g_2 +\frac{1}{2} + z)}
           {\text{sh} \frac{\alpha}{2}(\frac{1}{2} + z)}
      \frac{\text{ch} \frac{\alpha}{2} (g_3 +\frac{1}{2} + z)}
           {\text{ch} \frac{\alpha}{2}(\frac{1}{2} + z)}.
\end{eqnarray*}
The eigenvalues are determined by
\begin{displaymath}
E_r\: (y) = 2^r\!\!\!\sum_{\stackrel{J\subset \{
 1,\ldots ,n\} }{0\leq |J|\leq r}} \!\! (-1)^{r-|J|} \Bigl(
 \prod_{j\in J} \text{ch}(\alpha y_j) \!\!\!\! \sum_{r\leq
 l_1\leq\cdots\leq l_{r-|J|}\leq n}\!\!\!\!  \text{ch}(\alpha
 \rho_{l_1})\cdots \text{ch}(\alpha \rho_{l_{r-|J|}}) \Bigr) ,
\end{displaymath}
with
\begin{displaymath}
\rho_j = (n-j)g +(g_0+\cdots+g_3)/2
\end{displaymath}
($\rho =\sum_{j=1}^n \rho_j\, e_j$).
In the above formulas $|J|$ represents the number of elements of
$J\subset \{ 1,\ldots ,n \}$, and we have used the conventions that
empty products are equal to one,
$U_{K,\, p}(x)\equiv 1$ if $p=0$,
and that the second sum in $E_r(y)$ is equal to
one if $|J|=r$.

For $r=1$ the Difference equation~(\ref{diffeqr}) reduces to
a difference equation already considered by
Koornwinder \cite{koo:askey} and (for special parameters)
Macdonald \cite{mac:orthogonal}.

\section{Recurrence relations}\label{sec4}
To describe the recurrence relations for the multivariable
Askey-Wilson polynomials it is convenient to introduce dual parameters
$\hat{g}, \hat{g}_r$, which are related to the parameters
$g, g_r$ by
\begin{equation}\label{repa}
\hat{g} =g,
\;\;\;\;\;\;\;\;\;\;
\left(
\begin{array}{c} \hat{g}_0 \\ \hat{g}_1 \\ \hat{g}_2 \\ \hat{g}_3 \end{array}
                                                                   \right)
= \frac{1}{2} \left(
\begin{array}{rrrr}
1  &  1  &  1  &  1  \\
1  &  1  &  -1  &  -1  \\
1  &  -1  &  1  &  -1  \\
1  &  -1  &  -1  &  1
\end{array}             \right)
\left( \begin{array}{c} g_0 \\ g_1 \\ g_2 \\ g_3 \end{array} \right) .
\end{equation}
Furthermore, instead of working with monic polynomials we will
pass to a different normalization by introducing
\begin{equation}
P_\lambda (x) = c_\lambda^{-1} p_\lambda (x) ,\;\;\;\;\;\;\;\;
c_\lambda = \hat{\Delta}^+(\rho +\lambda)/\hat{\Delta}^+(\rho),
\end{equation}
where $\hat{\Delta}^+(x)$ is defined as in
Eq.~(\ref{weightp})
but with
the parameters $g, g_r$ replaced by the dual parameters
$\hat{g}, \hat{g}_r$. Similarly, we will also use
the notation $\hat{p}_\lambda (x)$ and $\hat{P}_\lambda (x)$
($=\hat{c}_\lambda^{-1}\hat{p}_\lambda (x)$ with
$\hat{c}_\lambda =\Delta^+(\hat{\rho}+\lambda )/\Delta^+(\hat{\rho})$)
for the corresponding dual polynomials (in which
$g, g_r$ again gets replaced by $\hat{g}, \hat{g}_r$).

It was conjectured by Macdonald \cite{mac:some} that the renormalized
multivariable Askey-Wilson polynomials $P_\lambda (x)$ satisfy the relation
\begin{equation}\label{dualr}
P_\lambda (\hat{\rho} +\mu)= \hat{P}_\mu (\rho
+\lambda),\;\;\;\;\;\;\;\;\;\;
\lambda ,\mu \in \Lambda.
\end{equation}
Recently, this conjecture has been proved by
Cherednik \cite{che:macdonald's}
for special parameters (corresponding to reduced root systems)
and in \cite{die:self-dual} for more general
parameters subject to the self-duality condition
\begin{equation}\label{self-dual}
g_0-g_1-g_2-g_3=0
\end{equation}
(implying $\hat{g}_r = g_r$).

The crucial point is now that the Difference equations~(\ref{diffeqr})
combined with Relation (\ref{dualr}) give rise to a
system of recurrence relations for the renormalized multivariable
Askey-Wilson polynomials. These recurrence relations read explicitly
\begin{equation}\label{recr}
\hat{E}_r (x)\: P_{\lambda} (x)=\!\!\!\!\!\!\!
\sum_{\stackrel{J\subset \{ 1,\ldots ,n\} ,\, 0\leq|J|\leq r}
               {\varepsilon_j=\pm 1,\; j\in J;\;
                \lambda +e_{\varepsilon J} \in \Lambda}} \!\!\!\!
\hat{U}_{J^c,\, r-|J|}(\rho +\lambda)\,
\hat{V}_{\varepsilon J,\, J^c}(\rho +\lambda)\,
P_{\lambda +e_{\varepsilon J}} (x) ,
\end{equation}
$r=1,\ldots ,n$ (where the superscript hats have again been employed to
indicate that parameters have been replaced by dual parameters).

To arrive at Eq.~(\ref{recr}) one first substitutes $x=\rho +\lambda$
in the difference equation
$\hat{D}_r \hat{P}_\mu=\hat{E}_r(\hat{\rho}+\mu) \hat{P}_\mu $
for the dual polynomials $\hat{P}_\mu$. Invoking of
Relation~(\ref{dualr}) and using that
$\hat{V}_{\varepsilon J,\, J^c}(\rho +\lambda)=0$ if
$\lambda +e_{\varepsilon J} \not\in \Lambda$
then leads to Eq.~(\ref{recr}) for
$x = \hat{\rho}+\mu$ with $\mu \in \Lambda$.
However, since Eq.~(\ref{recr})
describes a relation between trigonometric polynomials, knowing
that the relation is satisfied for all $\hat{\rho}+\mu$, $\mu \in
\Lambda$, is in fact sufficient to conclude that equality must hold
identically for all values of $x$.

\section{Orthonormalization}\label{sec5}
The difference operators $D_r$ (\ref{ados}) are symmetric with
respect to the inner product $\langle\cdot ,\cdot\rangle_\Delta$
(\ref{ip}) and the functions $E_r(y)$, determining the eigenvalues
in Eq.~(\ref{diffeqr}), separate the points
$\rho +\lambda$, $\lambda \in \Lambda$
\cite{die:commuting}. Hence, it follows that the eigenfunctions
$p_\lambda (x)$ are orthogonal with respect to
$\langle\cdot ,\cdot\rangle_\Delta$ (a priori the definition of the
polynomials guarantees only that
$\langle p_\mu ,p_\lambda \rangle_\Delta =0$ if $\mu < \lambda$).

It turns out that the Recurrence relations (\ref{recr}) may be employed to
compute the normalization constants turning the polynomials into an
orthonormal system. More precisely, if one works out both sides of the
identity
\begin{displaymath}
\langle \hat{E}_r P_\lambda , P_{\lambda +\omega_r}\rangle_\Delta =
\langle P_\lambda , \hat{E}_r P_{\lambda +\omega_r}\rangle_\Delta ,
\;\;\;\;\;\;\;\; \omega_r \equiv e_1+\cdots +e_r ,
\end{displaymath}
using Recurrence relation (\ref{recr}) and the orthogonality of the
polynomials, then one arrives at a relation between
$\langle P_\lambda , P_{\lambda}\rangle_\Delta$ and
$\langle P_{\lambda +\omega_r} ,
P_{\lambda +\omega_r}\rangle_\Delta $.
By iterating this relation one can express
$\langle P_\lambda , P_{\lambda}\rangle_\Delta$ in terms of
$\langle 1 , 1\rangle_\Delta$ (which corresponds to $\lambda =0$).
Since the value of $\langle 1 , 1\rangle_\Delta$ is known from the
work of Gustafson \cite{gus:generalization} (see also
\cite{kad:proof}), this solves the question of determining the
orthonormalization constants.
The answer, finally, reads \cite{die:self-dual}:
\begin{equation}\label{eval}
\langle p_\lambda , p_{\lambda}\rangle_\Delta =
2^n n! \:\hat{\Delta}^+ (\rho +\lambda ) \hat{\Delta}^- (\rho +\lambda) ,
\end{equation}
with $\hat{\Delta}^+(x)$ taken the same as in Section~\ref{sec4} and
\begin{equation}
\hat{\Delta}^-(x) = \prod_{1\leq j< k \leq n}
\hat{d}^-_v(x_j+x_k)\, \hat{d}^-_v(x_j-x_k)
\prod_{1\leq j\leq n} \hat{d}^-_w(x_j),
\end{equation}
where
\begin{eqnarray*}
\hat{d}^-_v(z)\; &=&\; q^{ \hat{g} z/2}
\frac{(q^{(z+1)};\, q)_\infty}
     {(q^{(-\hat{g}+z+1)};\, q)_\infty},\\[1ex]
\hat{d}^-_w(z)&=& q^{(\hat{g}_0+\cdots+\hat{g}_3)z/2}
\frac {(q^{(z+1)},\,
      -q^{(z+1)},\,
       q^{(1/2+z)},\,
      -q^{(1/2+z)};\, q)_\infty    }
     {(q^{(-\hat{g}_0+z+1)},\,
      -q^{(-\hat{g}_1+z+1)},\,
       q^{(-\hat{g}_2+1/2+z)},\,
      -q^{(-\hat{g}_3+1/2+z)};\, q)_\infty    }  .
\end{eqnarray*}

{\em Remarks: i.} The orthogonality of the multivariable Askey-Wilson
polynomials was first proved by Koornwinder \cite{koo:askey} using
only the difference operator $D_1$. This proof is based on
the continuity of
$\langle p_\lambda , p_{\mu}\rangle_\Delta$ in the parameters
and the fact that
for generic parameters $E_1(\rho +\lambda)\neq E_1(\rho +\mu)$
if $\lambda \neq \mu$.

{\em ii.} The value for
$\langle p_\lambda , p_{\lambda}\rangle_\Delta$ in (\ref{eval})
has first been conjectured
by Macdonald \cite{mac:some} and was then proved by Cherednik
for special parameters (corresponding to reduced root systems)
using the representation theory
of affine Hecke algebras \cite{che:double}.
Recently, Macdonald announced that
Cherednik's Hecke-algebraic techniques can be extended to a proof of
Formula~(\ref{eval}) valid for general
parameters \cite{mac:affine} (see also \cite{nou:macdonald}).

{\em iii.} Our derivation of the Recurrence relations (\ref{recr}) hinges on
Relation (\ref{dualr}). Since the proof of the latter relation
presented in \cite{die:self-dual}
only covers the self-dual case, to date our elementary proof of the
Recurrence relations (\ref{recr}) (and consequently
Formula (\ref{eval})) is not complete unless the Self-duality condition
(\ref{self-dual}) is satisfied.

\acknowledgments
The author would like to thank Prof. M. Noumi for
explaning his results in
\cite{nou:macdonald} and Prof. T. Oshima for the
kind hospitality at the University of Tokyo.

\end{document}